\documentclass[journal]{IEEEtran}

\usepackage{atbegshi}
\AtBeginDocument{\AtBeginShipoutNext{\AtBeginShipoutDiscard}}
\addtocounter{page}{-1}
\usepackage{cite}
\usepackage{amsmath,amsfonts}
\usepackage{algorithm}
\usepackage{algorithmic}

\usepackage{array}
\ifCLASSOPTIONcompsoc
  \usepackage[caption=false,font=normalsize,labelfont=sf,textfont=sf]{subfig}
\else
  \usepackage[caption=false,font=footnotesize]{subfig}
\fi

\usepackage{hyperref}       
\usepackage{url}            
\usepackage{booktabs}       
\usepackage{amsfonts}       
\usepackage{nicefrac}       
\usepackage{microtype}      
\usepackage{multirow}
\usepackage{lipsum}
\usepackage{graphicx}
\usepackage{placeins}       
\usepackage{tabularx}       
\usepackage[nohyperlinks,printonlyused]{acronym}
\usepackage{lineno}
\usepackage{svg}
\usepackage{import}
\usepackage{xifthen}
\usepackage{pdfpages}
\usepackage{transparent}
\usepackage{color}
\usepackage{nomencl}
\usepackage{float}
\usepackage{etoolbox}
\usepackage{textcomp}
\usepackage{stfloats}
\usepackage{verbatim}

\newcolumntype{L}[1]{>{\raggedright\let\newline\\\arraybackslash\hspace{0pt}}m{#1}}
\newcolumntype{C}[1]{>{\centering\let\newline\\\arraybackslash\hspace{0pt}}m{#1}}
\newcolumntype{R}[1]{>{\raggedleft\let\newline\\\arraybackslash\hspace{0pt}}m{#1}}

\newcommand{\timeIndex}{t}
\newcommand{\timeBlock}{tb}
\newcommand{\tbcur}{tb^\text{sel}}
\newcommand{\linkIndex}{g}
\newcommand{\regionIndex}{r}
\newcommand{\techIndex}{i}
\newcommand{\scenIndex}{\omega}
\newcommand{\scencur}{\omega^\text{sel}}
\newcommand{\iterIndex}{k}
\newcommand{\iterBest}{\hat{k}}
\newcommand{\timeAll}{T}
\newcommand{\timeBlockAll}{TB}
\newcommand{\linkAll}{G}
\newcommand{\regionAll}{R}
\newcommand{\techAll}{I}
\newcommand{\techConvAll}{\techAll_\text{C}}
\newcommand{\techStorAll}{\techAll_\text{S}}
\newcommand{\techTransAll}{\techAll_\text{T}}
\newcommand{\scenAll}{\Omega}
\newcommand{\iterAll}{K}
\newcommand{\techStorage}{\techIndex_\text{S}}
\newcommand{\techConverter}{\techIndex_\text{C}}
\newcommand{\techConvStorage}{\techIndex_\text{S}}
\newcommand{\techTransport}{\techIndex_\text{T}}

\newcommand{\IndSys}{Z^\text{total}}
\newcommand{\Zexp}{Z^\text{exp}}
\newcommand{\Zexpstab}{Z^\text{exp,stab}}
\newcommand{\Zop}{Z^\text{op}_{\scenIndex}}
\newcommand{\Zlow}{Z^\text{lower}}
\newcommand{\Zup}{Z^\text{upper}}
\newcommand{\Zupglob}{Z^\text{upper,glob}}
\newcommand{\Zmp}{Z^\text{MP}}
\newcommand{\Zmpunstab}{Z^\text{MP,unstab}}
\newcommand{\Zsp}{Z^\text{SP}_{\scenIndex}}
\newcommand{\Zspl}{z^\text{SP}_{\iterIndex,\scenIndex}}
\newcommand{\Zsptb}{Z^\text{SP}_{\scenIndex,\timeBlock}}
\newcommand{\Zsptbl}{z^\text{SP}_{\iterIndex,\scenIndex,\timeBlock}}
\newcommand{\Thet}{\theta_{\scenIndex}}
\newcommand{\Thetatb}{\theta_{\scenIndex,\timeBlock}}
\newcommand{\pistor}{\lambda^\text{stor}_{\scenIndex,\regionIndex,\techStorage}}
\newcommand{\pistortime}{\pi^\text{stor}_{\scenIndex,\timeIndex,\regionIndex,\techStorage}}
\newcommand{\pistortb}{\lambda^\text{stor,TS}_{\scenIndex,\timeBlock,\regionIndex,\techStorage}}
\newcommand{\pistortbl}{\lambda^\text{stor,TS}_{\iterIndex,\scenIndex,\timeBlock,\regionIndex,\techStorage}}
\newcommand{\piconv}{\lambda^\text{conv}_{\scenIndex,\regionIndex,\techConverter}}
\newcommand{\piconvtime}{\pi^\text{conv}_{\scenIndex,\timeIndex,\regionIndex,\techConverter}}
\newcommand{\piconvtb}{\lambda^\text{conv,TS}_{\scenIndex,\timeBlock,\regionIndex,\techConverter}}
\newcommand{\piconvtbl}{\lambda^\text{conv,TS}_{\iterIndex,\scenIndex,\timeBlock,\regionIndex,\techConverter}}
\newcommand{\pitransal}{\lambda^\text{trans,al}_{\scenIndex,\linkIndex,\techTransport}}
\newcommand{\pitransaltime}{\pi^\text{trans,al}_{\scenIndex,\timeIndex,\linkIndex,\techTransport}}
\newcommand{\pitransaltb}{\lambda^\text{trans,al,TS}_{\scenIndex,\timeBlock,\linkIndex,\techTransport}}
\newcommand{\pitransaltbl}{\lambda^\text{trans,al,TS}_{\iterIndex,\scenIndex,\timeBlock,\linkIndex,\techTransport}}
\newcommand{\pitransag}{\lambda^\text{trans,ag}_{\scenIndex,\linkIndex,\techTransport}}
\newcommand{\pitransagtime}{\pi^\text{trans,ag}_{\scenIndex,\timeIndex,\linkIndex,\techTransport}}
\newcommand{\pitransagtb}{\lambda^\text{trans,ag,TS}_{\scenIndex,\timeBlock,\linkIndex,\techTransport}}
\newcommand{\pitransagtbl}{\lambda^\text{trans,ag,TS}_{\iterIndex,\scenIndex,\timeBlock,\linkIndex,\techTransport}}
\newcommand{\pistorfix}{\pi^\text{storfix}_{\scenIndex,\timeBlock,\regionIndex,\techStorage}}
\newcommand{\pistorfixprev}{\pi^\text{storfix,prev}_{\scenIndex,\timeBlock,\regionIndex,\techStorage}}
\newcommand{\pistorl}{\lambda^\text{stor}_{\iterIndex,\scenIndex,\regionIndex,\techStorage}}
\newcommand{\piconvl}{\lambda^\text{conv}_{\iterIndex,\scenIndex,\regionIndex,\techConverter}}
\newcommand{\pitransall}{\lambda^\text{trans,al}_{\iterIndex,\scenIndex,\linkIndex,\techTransport}}
\newcommand{\pitransagl}{\lambda^\text{trans,ag}_{\iterIndex,\scenIndex,\linkIndex,\techTransport}}
\newcommand{\pistorfixl}{\lambda^\text{storfix}_{\iterIndex,\scenIndex,\timeBlock,\regionIndex,\techStorage}}

\newcommand{\pistorfixprevlnext}{\lambda^\text{storfix,prev}_{\iterIndex,\scenIndex,\timeBlock+1,\regionIndex,\techStorage}}
\newcommand{\pibdcut}{\pi^\text{cut}_{\iterIndex,\scenIndex,\timeBlock}}
\newcommand{\pibdcutother}{\pi^\text{cut}_{\iterIndex',\scenIndex,\timeBlock}}
\newcommand{\prob}{prob_{\scenIndex}}
\newcommand{\Cost}{m}
\newcommand{\CostInvConv}{\Cost^{\text{conv,inv}}_{\regionIndex,\techConverter}}
\newcommand{\CostFixConv}{\Cost^{\text{conv,fix}}_{\regionIndex,\techConverter}}
\newcommand{\CostInvStor}{\Cost^{\text{stor,inv}}_{\regionIndex,\techStorage}}
\newcommand{\CostFixStor}{\Cost^{\text{stor,fix}}_{\regionIndex,\techStorage}}
\newcommand{\CostInvTrans}{\Cost^{\text{trans,inv}}_{\linkIndex,\techTransport}}
\newcommand{\CostFixTrans}{\Cost^{\text{trans,fix}}_{\linkIndex,\techTransport}}
\newcommand{\CostVarConv}{\Cost^{\text{conv,var}}_{\techConverter}}

\newcommand{\CostFuel}{\Cost^{\text{fuel}}_{\techConverter}}

\newcommand{\CostUnmet}{\Cost^{\text{unserved}}}
\newcommand{\Capacity}{C}
\newcommand{\CapacityConvMin}{c^\text{conv,min}_{\regionIndex,\techConverter}}
\newcommand{\CapacityConvMax}{c^\text{conv,max}_{\regionIndex,\techConverter}}
\newcommand{\CapacityStorMin}{c^\text{stor,min}_{\regionIndex,\techStorage}}
\newcommand{\CapacityStorMax}{c^\text{stor,max}_{\regionIndex,\techStorage}}
\newcommand{\CapacityTransMin}{c^\text{trans,min}_{\linkIndex,\techTransport}}
\newcommand{\CapacityTransMax}{c^\text{trans,max}_{\linkIndex,\techTransport}}

\newcommand{\CapacityStor}{\Capacity^\text{stor}_{\regionIndex,\techStorage}}
\newcommand{\CapacityTrans}{\Capacity^\text{trans}_{\linkIndex,\techTransport}}
\newcommand{\CapacityConv}{\Capacity^\text{conv}_{\regionIndex,\techConverter}}
\newcommand{\CapacityStorL}{c^\text{stor}_{\iterIndex,\regionIndex,\techStorage}}
\newcommand{\CapacityTransL}{c^\text{trans}_{\iterIndex,\linkIndex,\techTransport}}
\newcommand{\CapacityConvL}{c^\text{conv}_{\iterIndex,\regionIndex,\techConverter}}
\newcommand{\CapacityStorLBest}{c^\text{stor}_{\iterBest,\regionIndex,\techStorage}}
\newcommand{\CapacityTransLBest}{c^\text{trans}_{\iterBest,\linkIndex,\techTransport}}
\newcommand{\CapacityConvLBest}{c^\text{conv}_{\iterBest,\regionIndex,\techConverter}}
\newcommand{\CapacityConvStor}{\Capacity^\text{conv}_{\regionIndex,\techConvStorage}}
\newcommand{\Dispatch}{P_{\scenIndex,\timeIndex,\regionIndex,\techConverter}}

\newcommand{\StorLevel}{L_{\scenIndex,\timeIndex,\regionIndex,\techStorage}}
\newcommand{\StorLevelPrev}{L_{\scenIndex,\timeIndex-1,\regionIndex,\techStorage}}

\newcommand{\Charge}{S^{\text{in}}_{\scenIndex,\timeIndex,\regionIndex,\techStorage}}
\newcommand{\Discharge}{S^{\text{out}}_{\scenIndex,\timeIndex,\regionIndex,\techStorage}}
\newcommand{\StorLoss}{S^{\text{loss}}_{\scenIndex,\timeIndex,\regionIndex,\techStorage}}

\newcommand{\Fuel}{J_{\scenIndex,\timeIndex,\regionIndex,\techConverter}}
\newcommand{\UnmetDemand}{D^\text{unserved}_{\scenIndex,\timeIndex,\regionIndex}}
\newcommand{\UnmetStorPlus}{L^\text{+}_{\scenIndex,\timeIndex,\regionIndex,\techStorage}}
\newcommand{\UnmetStorMinus}{L^\text{-}_{\scenIndex,\timeIndex,\regionIndex,\techStorage}}
\newcommand{\flow}{F^\text{al}_{\scenIndex,\timeIndex,\linkIndex,\techTransport}}
\newcommand{\flowagainst}{F^\text{ag}_{\scenIndex,\timeIndex,\linkIndex,\techTransport}}
\newcommand{\avail}{a_{\scenIndex,\timeIndex,\regionIndex,\techConverter}}
\newcommand{\demand}{d_{\scenIndex,\timeIndex,\regionIndex}}
\newcommand{\etaConv}{\eta_{\techConverter}}
\newcommand{\etaStorIn}{\eta^{\text{stor,in}}_{\techStorage}}
\newcommand{\etaStorOut}{\eta^{\text{stor,out}}_{\techStorage}}

\newcommand{\linkmatrix}{M^\text{line}_{\linkIndex,\regionIndex}}

\newcommand{\storagefixedVar}{L^\text{fix}_{\scenIndex,\timeBlock,\regionIndex,\techStorage}}
\newcommand{\storagefixedL}{l^\text{fix}_{\iterIndex,\scenIndex,\timeBlock,\regionIndex,\techStorage}}

\newcommand{\tbmatrix}{M^\text{last}_{\timeBlock,\timeIndex}}

\newcommand{\tballmatrix}{M^\text{time}_{\timeBlock,\timeIndex}}
\newcommand{\lbm}{\gamma}
\newcommand{\delt}{\delta_{\iterIndex,\scenIndex,\timeBlock}}
\newcommand{\deltother}{\delta_{\iterIndex',\scenIndex,\timeBlock}}
\newcommand{\act}{active_{\iterIndex,\scenIndex,\timeBlock}}
\newcommand{\actother}{active_{\iterIndex',\scenIndex,\timeBlock}}
\newcommand{\epsconv}{\varepsilon^\text{converge}}
\newcommand{\epsactive}{\varepsilon^\text{active}}

\begin{document}

\title{Accelerating Stochastic Energy System Optimization Models: Temporally Split Benders Decomposition}

\author{Shima~Sasanpour,
        Manuel~Wetzel,
        Karl-Kiên~Cao,
        Hans~Christian~Gils,
        Andrés~Ramos}
\thanks{Shima~Sasanpour, Manuel~Wetzel, Karl-Kiên~Cao and Hans~Christian~Gils are with the German Aerospace Center (DLR), Institute of Networked Energy Systems, Stuttgart, Germany. Andrés~Ramos is with the Institute for Research in Technology (IIT),
School of Engineering (ICAI), Universidad Pontificia Comillas, Madrid,
Spain. Corresponding author: \texttt{Shima.Sasanpour@dlr.de}} \\

\maketitle

\begin{abstract}
Stochastic programming can be applied to consider uncertainties in energy system optimization models for capacity expansion planning. However, these models become increasingly large and time-consuming to solve, even without considering uncertainties. For two-stage stochastic capacity expansion planning problems, Benders decomposition is often applied to ensure that the problem remains solvable. Since stochastic scenarios can be optimized independently within subproblems, their optimization can be parallelized. However, hourly-resolved capacity expansion planning problems typically have a larger temporal than scenario cardinality. Therefore, we present a temporally split Benders decomposition that further exploits the parallelization potential of stochastic expansion planning problems. A compact reformulation of the storage level constraint into linking variables ensures that long-term storage operation can still be optimized despite the temporal decomposition. We demonstrate this novel approach with model instances of the German power system with up to 87 million rows and columns. Our results show a reduction in computing times of up to 60\% and reduced memory requirements. Additional enhancement strategies and the use of distributed memory on high-performance computers further improve the computing time by over 80\%.
\end{abstract}



\renewcommand\nomgroup[1]{%
  \item[\bfseries
  \ifstrequal{#1}{D}{Sets}{%
  \ifstrequal{#1}{P}{Parameters}{%
  \ifstrequal{#1}{O}{Dual variables and subgradients}{%
  \ifstrequal{#1}{L}{First-stage variables}{\ifstrequal{#1}{M}{Second-stage variables}{}}}}}%
]}

\makenomenclature
\nomenclature[D,11]{$\timeIndex \in \timeAll$}{Time step index}
\nomenclature[D,08]{$\regionIndex \in \regionAll$}{Region index}
\nomenclature[D,01]{$\linkIndex \in \linkAll$}{Grid line index}
\nomenclature[D,02]{$\techIndex \in \techAll$}{Technology index}
\nomenclature[D,06]{$\iterIndex \in \iterAll$}{Iteration index}
\nomenclature[D,07]{$\iterBest$}{Iteration of current best solution}
\nomenclature[D,04]{$\techStorage \in \techStorAll$}{Storage technology index}
\nomenclature[D,03]{$\techConverter \in \techConvAll$}{Converter technology index}
\nomenclature[D,05]{$\techTransport \in \techTransAll$}{Transmission technology index}
\nomenclature[D,12]{$\timeBlock \in \timeBlockAll$}{Time block index}
\nomenclature[D,09]{$\scenIndex \in \scenAll$}{Scenario index}

\nomenclature[L,36]{$\IndSys$}{Total system cost}
\nomenclature[L,04]{$\CapacityTrans$}{Transmission capacity}
\nomenclature[L,01]{$C^\text{conv/stor}_{\regionIndex,\techIndex_\text{C/S}}$}{Converter/storage capacity}
\nomenclature[M,23]{$\Dispatch$}{Converter dispatch}
\nomenclature[M,10]{$\StorLevel$}{Storage level}
\nomenclature[M,11]{$\StorLevelPrev$}{Storage level of previous time step}
\nomenclature[M,25]{$S^{\text{in/out/loss}}_{\scenIndex,\timeIndex,\regionIndex,\techStorage}$}{Storage charging/discharging/loss}
\nomenclature[M,07]{$F^\text{al/ag}_{\scenIndex,\timeIndex,\linkIndex,\techTransport}$}{Power flow along/against line}
\nomenclature[M,09]{$\Fuel$}{Utilized fuel}
\nomenclature[L,28]{$\Zexp$}{Expansion cost}
\nomenclature[M,33]{$\Zop$}{Operational cost in scenario $\scenIndex$}
\nomenclature[L,30]{$Z^\text{lower/upper}$}{Lower/upper bound of system cost}
\nomenclature[L,38]{$\Zupglob$}{Global upper bound of system cost}
\nomenclature[L,31]{$\Zmp$}{Cost of MP}
\nomenclature[M,34]{$\Zsp$}{Cost of SP}
\nomenclature[M,35]{$\Zsptb$}{Cost of temporally split SP}
\nomenclature[L,39]{$\Thet$}{Approximation of SP cost}
\nomenclature[L,40]{$\Thetatb$}{Approximation of temporally split SP cost}
\nomenclature[M,05]{$\UnmetDemand$}{Unserved demand}
\nomenclature[M,12]{$L^\text{+/-}_{\scenIndex,\timeIndex,\regionIndex,\techStorage}$}{Slack adding/reducing storage level}
\nomenclature[L,22]{$\storagefixedVar$}{Fixed storage level}

\nomenclature[O,04]{$\pibdcut$}{Dual variable of optimality cut}
\nomenclature[O,16]{$\pistor$}{Subgradient of storage level constraint}
\nomenclature[O,05]{$\pistortime$}{Dual variable of storage level constraint}
\nomenclature[O,12]{$\piconv$}{Subgradient of power generation constraint}
\nomenclature[O,01]{$\piconvtime$}{Dual variable of power generation constraint}
\nomenclature[O,21]{$\pitransal$}{Subgradient of power transmission constraint along line}
\nomenclature[O,10]{$\pitransaltime$}{Dual variable of power transmission constraint along line}
\nomenclature[O,24]{$\pitransag$}{Subgradient of power transmission constraint against line}
\nomenclature[O,11]{$\pitransagtime$}{Dual variable of power transmission constraint against line}
\nomenclature[O,08]{$\pistorfix$}{Dual variable of fixed storage level constraint for last time step of time block}
\nomenclature[O,09]{$\pistorfixprev$}{Dual variable of fixed storage level constraint for last time step of previous time block}

\nomenclature[P,01]{$\avail$}{Power plant availability}
\nomenclature[P,02]{$\act$}{Parameter indicating if a cut is active}
\nomenclature[P,03]{$\CapacityConvL$}{Optimized converter capacity in iteration $\iterIndex$}
\nomenclature[P,04]{$c^\text{conv/stor,min/max}_{\regionIndex,\techIndex_\text{C/S}}$}{Min./max.\ converter/storage capacity}
\nomenclature[P,06]{$\CapacityStorL$}{Optimized storage capacity in iteration $\iterIndex$}
\nomenclature[P,09]{$\CapacityTransL$}{Optimized transmission capacity in iteration $\iterIndex$}
\nomenclature[P,11]{$c^\text{trans,min/max}_{\linkIndex,\techTransport}$}{Min./max.\ transmission capacity}
\nomenclature[P,12]{$\demand$}{Electricity demand}
\nomenclature[P,15]{$\Cost^{\text{conv/stor,inv}}_{\regionIndex,\techIndex_\text{C/S}}$}{Specific converter/storage investment cost}
\nomenclature[P,16]{$\Cost^{\text{conv/stor,fix}}_{\regionIndex,\techIndex_\text{C/S}}$}{Specific converter/storage fixed O\&M cost}
\nomenclature[P,17]{$\CostVarConv$}{Specific variable O\&M cost}
\nomenclature[P,19]{$\CostFuel$}{Specific fuel cost}
\nomenclature[P,22]{$\CostInvTrans$}{Investment cost of transmission technology $\techTransport$}
\nomenclature[P,23]{$\CostFixTrans$}{Fixed O\&M cost of transmission technology $\techTransport$}
\nomenclature[P,24]{$\CostUnmet$}{Penalty cost}
\nomenclature[P,25]{$\linkmatrix$}{Matrix linking transmission lines and regions}
\nomenclature[P,26]{$\tbmatrix$}{Matrix linking time blocks to their last time step}
\nomenclature[P,27]{$\tballmatrix$}{Matrix linking time blocks and time steps}
\nomenclature[P,28]{$\prob$}{Probability of scenario $\scenIndex$}
\nomenclature[P,31]{$\Zspl$}{Cost of SP in iteration $\iterIndex$}
\nomenclature[P,32]{$\Zsptbl$}{Cost of temporally split SP in iteration $\iterIndex$}
\nomenclature[P,33]{$\beta$}{Level weighting parameter}
\nomenclature[P,34]{$\delt$}{Number of the iteration when the cut was created or lastly binding}
\nomenclature[P,35]{$\epsactive$}{Threshold indicating if a cut is binding}
\nomenclature[P,36]{$\epsconv$}{Convergence tolerance}
\nomenclature[P,37]{$\etaConv$}{Converter efficiency}
\nomenclature[P,38]{$\eta^{\text{stor,in/out}}$}{Charging/Discharging efficiency}
\nomenclature[P,40]{$\phi$}{Number of iterations a cut needs to be unbinding to be deactivated}
\nomenclature[P,41]{$\lbm$}{Level parameter}

\printnomenclature[2cm]

\IEEEpeerreviewmaketitle

\section{Introduction}
\label{sec:intro}
\IEEEPARstart{R}{educing} the environmental and climate damage caused by the global energy supply is a key challenge of our time. Energy system optimization models (ESOMs) for capacity expansion planning (CEP) can support decision makers and system planners to determine the least-cost decarbonized energy systems. Typically, these models optimize the installed capacities of the technologies used in the system and their operation over the course of a year by minimizing the total cost of the system. However, they are based on various assumptions about future developments, such as weather-based power feed-in from variable renewable energy technologies and technology costs, for which the exact values are uncertain.

These uncertainties are often ignored in CEP, even though they can have a significant impact on the optimized infrastructure of the energy system. Yue et al.\ identify four potential ways to systematically consider uncertainties in ESOMs \cite{yue2018review}. Monte Carlo analysis (MCA) has the ability to analyze a large scenario space covering different uncertainties. However, the large range of possible results makes it difficult to derive concrete advice from them.
Modeling to generate alternatives (MGA) enables the consideration of energy systems of slightly higher cost but highly different infrastructures by maximizing the distance to the cost-optimal solution. Similarly to MCA, the wide range of possible results can show a lot of options, but not one optimized energy system with the desired properties to consider different types of uncertainties.
Robust optimization can range from a worst-case optimization to a risk-averse optimization where a budget of uncertainty is defined. Due to this risk-aversion, pessimistic scenarios with low probability will still largely influence the optimization, making the energy system more expensive.
Finally, stochastic programming has the advantage that, on the one hand, it results in one optimization strategy for all considered uncertainties, similar to robust optimization. On the other hand, all stochastic scenarios are assigned a probability, which allows the modeler to assign lower probabilities to less likely scenarios while still being able to consider them in the optimization. This allows us to determine an energy system that can hedge the risk of the considered uncertainties, which is why we aim to consider uncertainties by applying stochastic programming in our optimization.

The informative value of the analyses with ESOMs increases with their ability to map the details and scope of the real system and its operation. This drives the ambition to consider many technologies and energy carriers in high detail, to consider a large spatial granularity for the representation of energy networks, and also to optimize the use of these technologies for each of the 8760 hours of the year. In models for the consideration of large-scale systems, i.e.\ national or continental, this inevitably leads to very large optimization problems that are complex and time-consuming to solve using standard methods. This is exacerbated by efforts to consider uncertainties in CEP through stochastic programming. In this study we consider large-scale optimization problems with up to 87 million rows and columns.

Benders decomposition (BD) is an established method to split the problem into smaller parts that are solved iteratively until the optimal solution is found \cite{benders1962partitioning}. The master problem (MP) typically optimizes the complicating variables such as linking variables representing the decisions on the expansion of technologies in the system, e.g.\ power plants, storage technologies, or networks. When stochastic programming problems are considered, the uncertainties related to the operation of the energy system, e.g.\ electricity demand \cite{soares2017two} or the availability of variable renewable resources \cite{goke2024stabilized}, can be accounted for in the subproblems (SPs). Here, each SP represents one stochastic scenario. This has the advantage that the SPs are independent of each other and can therefore be solved in parallel. However, the BD algorithm in its classic formulation may still be very slow, since the convergence of the problem can take up a lot of iterations. Rahmaniani et al.\ summarize different approaches to improve the solving speed and convergence of the BD, called enhancement strategies \cite{rahmaniani2017benders}. As stated by Göke et al.\, many enhancement strategies are related to the improved calculation of the MP, since this is the complicating part in many optimization problems \cite{goke2024stabilized}. Crainic et. al. adjust the decomposition strategy of the considered two-stage stochastic problem by including explicit information from the SPs within the MP, improving the efficiency and stability of the solution process despite the increased difficulty of the MP \cite{crainic2016partial}. Rahmaniani et al.\ describe the modification of the decomposition strategy as a promising enhancement approach, although research in that field is rather limited \cite{rahmaniani2017benders}.

In CEP, the dispatch problems that are solved within the SPs tend to get comparatively large due to the many time steps that are considered. This results in rather computationally expensive SPs as compared to the MPs. Therefore, we propose a decomposition of the SPs not only along the scenario set dimension but additionally along the time dimension. This could, on the one hand, decrease the size of the SPs, making them easier and faster to solve. On the other hand, the parallelization potential could be further exploited.
BD has been applied for CEP before. Grübler et al.\ performed a spatial and temporal decomposition on a deterministic CEP problem using BD, resulting in a reduction of the solving time compared to BD without time decomposition \cite{grubler2024applying}. Jacobson et al.\ introduce a BD approach with temporal decomposition for a deterministic mixed-integer linear programming CEP problem with an annual emission limit constraint. In order to still achieve the optimal solution, budgeting variables are used within the MP \cite{jacobson2024computationally}. However, neither studies consider storage technologies within their temporal decomposition. As interconnected future energy systems heavily rely on variable renewable energy sources, flexibility options such as storage technologies will become more relevant. These storage technologies result in linking constraints since the storage levels connect consecutive time steps. A temporal decomposition can therefore not optimize multi-day storage technologies, such as pumped hydro storage. While Pecci et al.\ introduce a first approach on the consideration of the multi-day storage technologies when temporally decomposing the SPs, their main focus is on the analysis of bundle methods for mixed-integer problems \cite{pecci2025regularized}.
Nested Benders decomposition is another approach that enables the temporal decomposition in CEP. However, due to the temporal hierarchy of the SPs the parallelization potential is usually rather limited \cite{santos2016new}.

This paper presents an improved BD algorithm for stochastic, very-large-scale CEP problems that applies time decomposition to reduce the solving time. A novel compact formulation for the additional optimization of the storage level of the last time step of each time block within the MP ensures that long-term storage technologies can be optimized despite the temporal decomposition. Therefore, the same optimal long-term storage operation can be achieved. This additional decomposition of all time steps into several time blocks facilitates the parallelization of smaller SPs. By integrating MPI (Message Passing Interface), the computation on distributed memory architectures becomes possible, enabling access to high-performance computing (HPC). Further enhancement strategies, such as the utilization of bundle methods, are integrated, resulting in a stabilized convergence of the algorithm. This results in considerable time savings compared to solving the deterministic equivalent (DEQ).

\section{Method}
\label{sec:Method}
\subsection{Energy system optimization framework REMix}
\label{subsec:REMix}
\noindent REMix is a GAMS-based open-source framework for optimizing the design and operation of energy systems \cite{wetzel2024remix}. The scope of the models built with REMix is very flexible and can include various energy carriers, such as power, heat, and synthetic fuels \cite{wetzel2023green}, and a geographical resolution ranging from country- \cite{sasanpour2021strategic} to transformer-substation level \cite{frey2024tackling}. In its basic form, REMix performs a deterministic linear optimization of one target year with hourly resolution for generation, storage, and transmission capacities. However, further advanced features are available, e.g.\ unit commitment and multi-year optimization with perfect foresight. 
Additionally, stochastic programming can be applied to generate and solve the DEQ, e.g.\ to model the uncertain power generation of variable renewable energy sources \cite{sasanpour2022quantifying}.

The objective function of the DEQ
\begin{align}
\label{eq:syscost}
    \IndSys & = \min \Zexp + \sum_{\scenIndex} \prob \Zop
\end{align}
minimizes the total system costs $\IndSys$, consisting of the cost for expansion $\Zexp$ and the expected cost for the operation $\Zop$ of the energy system, taking the probability $\prob$ of each scenario $\scenIndex$ into account. The expansion cost
\begin{align}
\label{eq:expcost}
    \Zexp & = \sum_{\regionIndex,\techConverter} (\CostInvConv + \CostFixConv) \CapacityConv \\ & \nonumber + \sum_{\regionIndex,\techStorage} (\CostInvStor + \CostFixStor) \CapacityStor \\ & \nonumber + \sum_{\linkIndex,\techTransport} (\CostInvTrans + \CostFixTrans) \CapacityTrans
\end{align}
represents the annuity and fixed operation and maintenance (O\&M) costs for expanded converter, storage, and transmission technologies, while the operational cost in each scenario
\begin{align}
\label{eq:opcost}
    \Zop & = \sum_{\timeIndex,\regionIndex,\techConverter} \CostVarConv \Dispatch + \sum_{\timeIndex,\regionIndex,\techConverter} \CostFuel \Fuel \\ & \nonumber + \sum_{\timeIndex,\regionIndex} \CostUnmet \UnmetDemand, \forall \text{ } \scenIndex \in \scenAll
\end{align}
includes variable O\&M cost for power generation, and fuel cost.
If the demand can not be met, this is accounted for by additional penalty costs.

The model is subject to a set of constraints. The capacities of converter, storage and transmission technologies $\CapacityConv$, $\CapacityStor$ and $\CapacityTrans$ are restricted by both lower and upper limits
\begin{align}
\label{eq:capacityconvmax}
   \CapacityConvMin & \leq \CapacityConv \leq \CapacityConvMax, \forall \text{ } \regionIndex \in \regionAll,\techConverter \in \techConvAll, \\
\label{eq:capacitystormax}
    \CapacityStorMin & \leq \CapacityStor \leq \CapacityStorMax, \forall \text{ } \regionIndex \in \regionAll,\techStorage \in \techStorAll, \\
\label{eq:capacitytransmax}
    \CapacityTransMin & \leq \CapacityTrans \leq \CapacityTransMax, \forall \text{ } \linkIndex \in \linkAll,\techTransport \in \techTransAll.
\end{align}
The lower limits can e.g.\ represent capacities that have been built in previous years and where the technical lifetime has not been exceeded, yet.
The expanded capacities
\begin{align}
\label{eq:dispatch}
    \CapacityConv \avail \geq & \text{ } \Dispatch : \piconvtime, \forall \text{ } \scenIndex \in \scenAll, \\ & \nonumber \timeIndex \in \timeAll,\regionIndex \in \regionAll,\techConverter \in \techConvAll, \\
\label{eq:storcap}
    \CapacityStor \geq & \text{ } \StorLevel : \pistortime, \forall \text{ } \scenIndex \in \scenAll, \\ & \nonumber \timeIndex \in \timeAll, \regionIndex \in \regionAll,\techStorage \in \techStorAll, \\
\label{eq:transpcap}
    \CapacityTrans \geq & \text{ } \flow : \pitransaltime, \forall \text{ } \scenIndex \in \scenAll, \\ & \nonumber \timeIndex \in \timeAll, \linkIndex \in \linkAll,\techTransport \in \techTransAll, \\
\label{eq:transpcapagainst}
    \CapacityTrans \geq & \text{ } \flowagainst : \pitransagtime, \forall \text{ } \scenIndex \in \scenAll, \\ & \nonumber \timeIndex \in \timeAll, \linkIndex \in \linkAll,\techTransport \in \techTransAll
\end{align}
limit the power generation of the power plants $\Dispatch$ (Eq. (\ref{eq:dispatch})), the storage level $\StorLevel$ (Eq. (\ref{eq:storcap})) and the power transmission along $\flow$ (Eq. (\ref{eq:transpcap})) and against $\flowagainst$ (Eq. (\ref{eq:transpcapagainst})) each line $\linkIndex$ from region $\regionIndex$ to region $\regionIndex'$, where $\piconvtime$, $\pistortime$, $\pitransaltime$ and $\pitransagtime$ represent the respective dual variables of the equations. When power plants are dispatched, their availability $\avail$ is taken into account, representing planned and unplanned unavailabilities in the case of conventional power plants and weather fluctuations in the case of renewable technologies. The availability of storage and transmission technologies can also be restricted. However, this is not considered in this study.

The storage level
\begin{align}
\label{eq:storlevel}
    \StorLevel & = \StorLevelPrev + \Charge \etaStorIn - \frac{\Discharge}{\etaStorOut} \\ & \nonumber -\StorLoss, \forall \text{ } \scenIndex \in \scenAll,\timeIndex \in \timeAll,\regionIndex \in \regionAll,\techStorage \in \techStorAll
\end{align}
depends on the storage level from the previous time step $\StorLevelPrev$, the amount of power that is charged $\Charge$ and discharged $\Discharge$ (taking charging and discharging efficiencies $\etaStorIn$ and $\etaStorOut$ into account) and storage losses $\StorLoss$. Due to a circular formulation of Eq. (\ref{eq:storlevel}), the storage level of the last time step is connected to the storage level of the first time step, so that the storage is not necessarily completely emptied at the end of the year.

Converter technologies do not only represent power plants, but they can also be utilized to charge and discharge storage technologies. E.g.\, in the case of a pumped hydro storage, the storage capacity $\CapacityStor$ determines the storage energy that can be stored. The converter capacity $\CapacityConvStor$ indicates how fast the storage can be charged (by pumps) and discharged (by turbines). Therefore, the amount that can be charged and discharged per time step 
\begin{align}
\label{eq:chargelimit}
    \Charge \leq \CapacityConvStor, \forall & \text{ } \scenIndex \in \scenAll,\timeIndex \in \timeAll, \\ & \nonumber \regionIndex \in \regionAll,\techConvStorage \in \techStorAll \\
\label{eq:dischargelimit}
    \Discharge \leq \CapacityConvStor, \forall & \text{ } \scenIndex \in \scenAll,\timeIndex \in \timeAll, \\ & \nonumber \regionIndex \in \regionAll,\techConvStorage \in \techStorAll
\end{align}
depends on the converter capacity of the storage technologies. The dispatch of the conventional power plants leads to the consumption of fuels 
\begin{align}
\label{eq:fuel}
    \Fuel = \frac{\Dispatch}{\etaConv}, \forall \text{ } \scenIndex \in \scenAll, & \timeIndex \in \timeAll,\regionIndex \in \regionAll,
\end{align}
taking the power plant efficiency $\etaConv$ into account. In our case, we only consider biomass and hydrogen-fueled power plants. Therefore, we assume that no carbon is emitted.

Finally, the demand
\begin{align}
\label{eq:balance}
     \demand & = \sum_{\techConverter} \Dispatch + \sum_{\techStorage} (\Discharge - \Charge \\ & \nonumber - \StorLoss) + \sum_{\linkIndex,\techTransport} \linkmatrix (\flowagainst - \flow) \\ & \nonumber + \UnmetDemand, \forall \scenIndex \in \scenAll, \timeIndex \in \timeAll,\regionIndex \in \regionAll
\end{align}
must be balanced with the supply in each time step $\timeIndex$ and region $\regionIndex$. The incidence matrix $\linkmatrix$ indicates which lines $\linkIndex$ start at region $\regionIndex$ with positive entries and which end in region $\regionIndex$ with negative entries.

\subsection{Benders decomposition for stochastic CEPs}
\label{subsec:Benders}
\begin{figure}[!t]
    \centering
    \includegraphics[width=90mm]{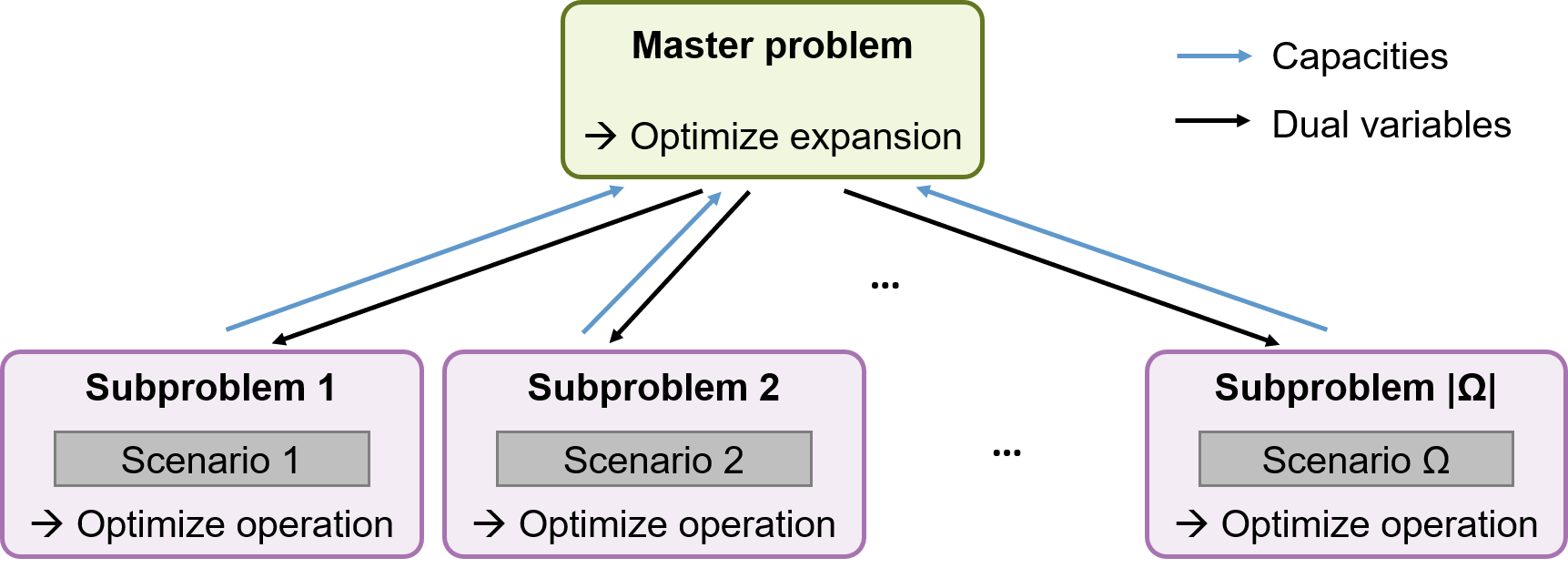}
    \caption{Benders decomposition for stochastic CEP problem. The MP optimizes the expansion of the energy system and communicates the optimized capacities to the SPs. Each SP optimizes the operation of one stochastic scenario. The SPs send dual variables back to the MP. Within the MP they are considered as optimality cuts to estimate the cost of the SPs. This process is repeated until the distance between the lower and upper bounds of the objective function is within the predefined tolerance.}
    \label{fig:Benders}
\end{figure}
\noindent The BD splits the stochastic CEP problem into one MP and several SPs. Fig. \ref{fig:Benders} shows the BD for a two-stage stochastic problem and the exchange of information between the MP and the SPs. The MP typically optimizes the linking variables. Within a stochastic CEP these are the expansion variables of the first stage, which limit the second-stage dispatch of the SPs (see Eq. (\ref{eq:dispatch}) - Eq. (\ref{eq:transpcapagainst})). The objective function of the MP
\begin{align}
\label{eq:objMP}
     \Zmp & = \min \Zexp + \sum_{\scenIndex} \Thet
\end{align}
minimizes the expansion cost from Eq. (\ref{eq:expcost}) and the estimated costs of the SPs. The cost of each SP (and therefore each stochastic scenario) is approximated as
\begin{align}
\label{eq:theta}
     \Thet & \geq \prob (\Zspl - \sum_{\regionIndex,\techStorage} (\CapacityStorL - \CapacityStor) \pistorl \\ & \nonumber - \sum_{\regionIndex,\techConverter} (\CapacityConvL - \CapacityConv) \piconvl \\ & \nonumber - \sum_{\linkIndex,\techTransport} (\CapacityTransL - \CapacityTrans) (\pitransall + \pitransagl)), \\ & \nonumber \forall \text{ } \iterIndex \in \iterAll, \scenIndex \in \scenAll.
\end{align}
The parameters $\CapacityStorL$, $\CapacityConvL$ and $\CapacityTransL$ store the optimized capacities of each passed iteration $\iterIndex$. The actual cost of SP $\scenIndex$ in iteration $\iterIndex$ is stored in $\Zspl$. Within the MP the multi-cut formulation is used. This means that each SP, and therefore stochastic scenario, generates one cut per iteration within the MP. This enhancement strategy is described in more detail in Section \ref{subsubsec:multi}. Additionally, Eq. (\ref{eq:capacityconvmax}) - Eq. (\ref{eq:capacitytransmax}) are taken into account as constraints for the MP.

Within each SP the capacity variables are fixed to the values optimized in the MP. Each SP minimizes the operation of one stochastic scenario
\begin{align}
\label{eq:objSP}
     \Zsp & = \min \Zop
\end{align}
subject to Eq. (\ref{eq:dispatch}) - Eq. (\ref{eq:balance}). After solving the SPs, the subgradients
\begin{align}
\label{eq:piconvsum}
    \piconv = & \sum_{\timeIndex} \piconvtime \avail, \forall \text{ } \scenIndex \in \scenAll, \\ & \nonumber \regionIndex \in \regionAll, \techConverter \in \techConvAll \\
\label{eq:pistorsum}
    \pistor = & \sum_{\timeIndex} \pistortime, \forall \text{ } \scenIndex \in \scenAll,\regionIndex \in \regionAll,\techStorage \in \techStorAll \\
\label{eq:pitransalsum}
    \pitransal = & \sum_{\timeIndex} \pitransaltime, \forall \text{ } \scenIndex \in \scenAll, \linkIndex \in \linkAll,\techTransport \in \techTransAll \\
\label{eq:pitransagsum}
    \pitransag = & \sum_{\timeIndex} \pitransagtime, \forall \text{ } \scenIndex \in \scenAll, \linkIndex \in \linkAll,\techTransport \in \techTransAll
\end{align}
are calculated by taking the sum over the time dimension of the dual variables related to the expansion variables (see Eq. (\ref{eq:dispatch}) - Eq. (\ref{eq:transpcapagainst})). This reduces the size of the parameters that need to be communicated to the MP. For the subgradient related to the power generation constraint, the availability factor $\avail$ is additionally considered. The subgradients are sent back to the MP and considered in Eq. (\ref{eq:theta}) in the next iteration. Since we consider unserved demand with penalty costs in Eq. (\ref{eq:opcost}), the SPs can not become infeasible. Therefore, only optimality cuts and no feasibility cuts are added to the MP. The algorithm is finished as soon as the distance between the lower bound $\Zlow$ and upper bound $\Zup$ of the objective function is within the convergence tolerance $\epsconv > 1 -\Zlow/\Zup$.

Since the SPs are independent of each other, they can be solved in parallel. This parallelization can help to reduce the solving time of the SPs, which are typically more time-consuming to solve than the MP \cite{goke2024stabilized}.

\subsection{Temporally split Benders decomposition}
\label{subsec:timesplitting}
\begin{figure*}[!ht]
    \centering
    \includegraphics[width=130mm]{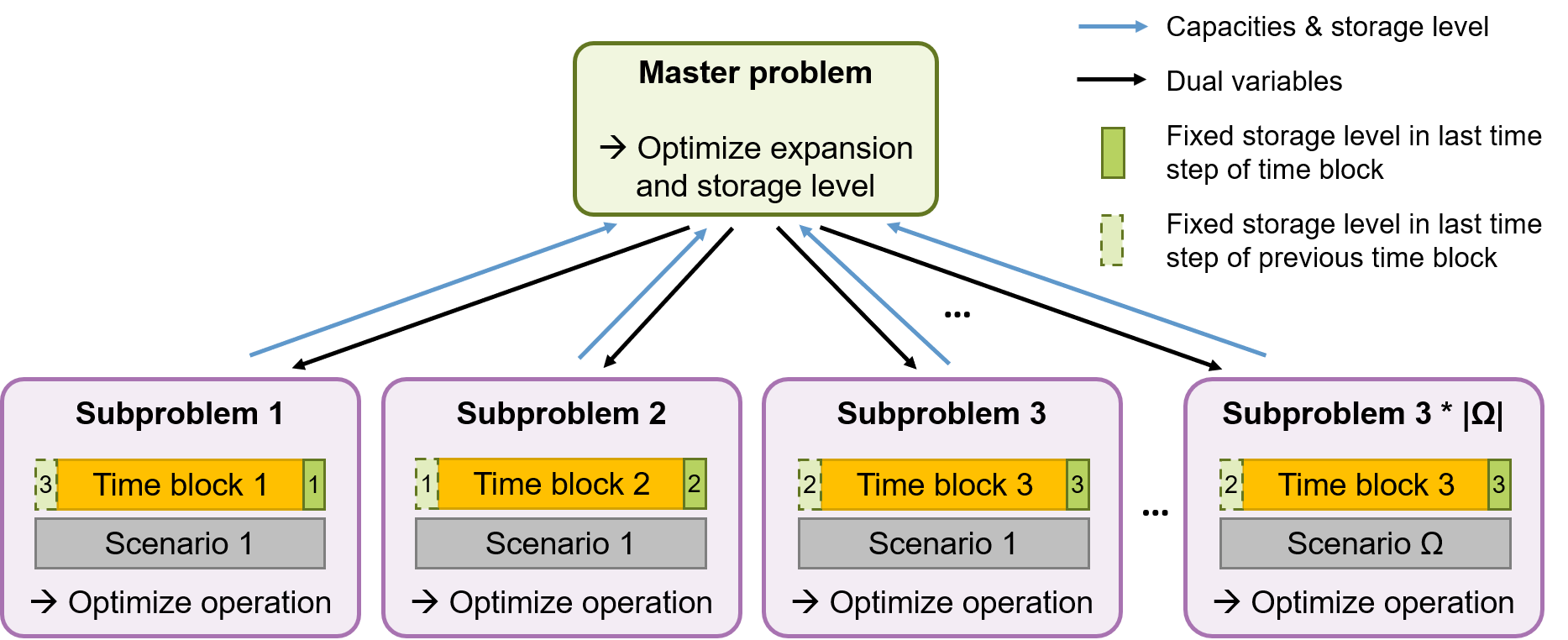}
    \caption{Temporally split Benders decomposition, exemplary for three time blocks per stochastic scenario. The MP optimizes the expansion and the storage level of the last time step of each time block and scenario. Each SP optimizes the operation of one time block of a scenario. The capacity and the storage level of the last time step of each time block are fixed. After solving, the SPs send dual variables related to the expansion and the storage level of the last time step of the respective and previous time block back to the MP.}
    \label{fig:timesplitting}
\end{figure*}
\noindent The high temporal resolution of up to 8760 time steps representing each hour of the year makes the time dimension usually much larger than the scenario dimension in ESOMs. To reduce the solving time of the SPs and to further increase the parallelization, we extend the decomposition strategy and divide the stochastic scenarios into several time blocks $\timeBlock$. However, the consideration of energy storage levels links all time steps of the model (see Eq. (\ref{eq:storlevel})). The SPs for the same scenario and different time blocks are thus not independent from each other. Therefore, we adjust the MP such that it not only optimizes the capacities of the energy system but also the storage level of the last time step of each time block, transforming the linking equations Eq. (\ref{eq:storlevel}) into $|\timeBlockAll|$ linking variables. Fig. \ref{fig:timesplitting} shows the exchange of information between the MP and the SPs of the temporally split Benders decomposition (TSBD) if each scenario is split into three time blocks. Each SP optimizes the operation of one time block ($\tbcur$) and one scenario ($\scencur$), resulting in $|\timeBlockAll| \cdot |\scenAll|$ SPs that can be solved in parallel. Besides the capacities, the storage level of the last time step of each time block
\begin{align}
\label{eq:storfix}
    \sum_{\timeIndex} &\tbmatrix \StorLevel = \storagefixedVar : \pistorfix , \\ & \nonumber \forall \text{ } \scenIndex \in \scencur, \timeBlock \in \tbcur, \regionIndex, \techStorage
\end{align}
is fixed. The matrix $\tbmatrix$ maps the last time step of a time block to the time block. The variable $\storagefixedVar$ represents the optimized storage level from the MP and is fixed in the SPs. Each SP additionally receives information on the storage level of the last time step of the previous time block
\begin{align}
\label{eq:storfixprev}
    \sum_{\timeIndex} &\tbmatrix \StorLevel = \storagefixedVar : \pistorfixprev, \\ & \nonumber \forall \text{ } \scenIndex \in \scencur, \timeBlock \in (\tbcur-1), \regionIndex, \techStorage
\end{align}
impacting the storage level of the first time step of the respective time block. The storage level of all but the last time step can be optimized within the SP, taking the fixed storage levels into account. Prior approaches had to optimize both the storage level of the first and last time steps within a time block to consider multi-day storage technologies in temporally split SPs \cite{pecci2025regularized}. The objective function of each SP
\begin{align}
\label{eq:objSPTS}
     \Zsptb & = \min \sum_{\timeIndex} \tballmatrix (\sum_{\regionIndex,\techConverter} \CostVarConv \Dispatch \\ & \nonumber+ \sum_{\regionIndex,\techConverter} \CostFuel \Fuel + \sum_{\regionIndex} \CostUnmet (\UnmetDemand \\ & \nonumber + \UnmetStorPlus + \UnmetStorMinus)), \forall \text{ } \scenIndex \in \scencur, \timeBlock \in \tbcur
\end{align}
minimizes the operational cost of one scenario $\scencur$ and one time block $\tbcur$. The matrix $\tballmatrix$ links all time steps to the respective time block. To receive the new subgradients
\begin{align}
\label{eq:piconvsumtb}
    \piconvtb = \sum_{\timeIndex} & \tballmatrix \piconvtime \avail, \\ & \nonumber \forall \text{ } \scenIndex \in \scenAll, \timeBlock \in \timeBlockAll, \regionIndex \in \regionAll,\techConverter \in \techConvAll \\
\label{eq:pistorsumtb}
    \pistortb = \sum_{\timeIndex} & \tballmatrix \pistortime, \forall \text{ } \scenIndex \in \scenAll, \\ & \nonumber \timeBlock \in \timeBlockAll, \regionIndex \in \regionAll,\techStorage \in \techStorAll \\
\label{eq:pitransalsumtb}
    \pitransaltb = \sum_{\timeIndex} & \tballmatrix \pitransaltime, \forall \text{ } \scenIndex \in \scenAll, \\ & \nonumber \timeBlock \in \timeBlockAll, \linkIndex \in \linkAll,\techTransport \in \techTransAll \\
\label{eq:pitransagsumtb}
    \pitransagtb = \sum_{\timeIndex} & \tballmatrix \pitransagtime, \forall \text{ } \scenIndex \in \scenAll,\\ & \nonumber \timeBlock \in \timeBlockAll, \linkIndex \in \linkAll,\techTransport \in \techTransAll
\end{align}
the dual variables related to the fixed capacities are summed up over all time steps in a time block and are then provided to the MP. Additionally, the dual variables $\pistorfix$ and $\pistorfixprev$ of Eq. (\ref{eq:storfix}) and Eq. (\ref{eq:storfixprev}) related to the fixed storage level of the last time step of the considered time block $\tbcur$ and of the previous time block $\tbcur$-1 are provided to the MP. Since the dual variables are now available for each time block and scenario, the multi-cut formulation can be extended. Each SP can provide an optimality cut to the MP, and therefore, the improved multi-cut formulation generates not only one optimality cut for each scenario but also for each scenario and time block combination. The approximation of the cost of the SPs within the MP is therefore redefined as
\begin{align}
\label{eq:thetanew}
     \Thetatb & \geq \prob \\ & \nonumber (\Zsptbl - \sum_{\regionIndex,\techStorage} (\CapacityStorL - \CapacityStor) \pistortbl \\ & \nonumber - \sum_{\regionIndex,\techConverter} (\CapacityConvL - \CapacityConv) \piconvtbl \\ & \nonumber - \sum_{\linkIndex,\techTransport} (\CapacityTransL - \CapacityTrans) (\pitransaltbl + \pitransagtbl) \\ & \nonumber - \sum_{\timeBlock,\regionIndex,\techStorage} (\storagefixedL - \storagefixedVar) \\ & \nonumber (\pistorfixl + \pistorfixprevlnext)) : \pibdcut, \forall \text{ } \iterIndex \in \iterAll, \\ & \nonumber \scenIndex \in \scenAll, \timeBlock \in \timeBlockAll, (\iterIndex,\scenIndex,\timeBlock) \in \act.
\end{align}
The parameter $\act$ indicates if an optimality cut is active. The number of cuts added per iteration increases further, resulting in more information being provided to the MP per iteration. This can result in faster convergence, however, at the cost of a faster-growing MP. The objective function of the MP is updated to
\begin{align}
\label{eq:objMPTS}
     \Zmp & = \min \Zexp + \sum_{\scenIndex,\timeBlock} \Thetatb.
\end{align}

In addition to the slack variable $\UnmetDemand$ for unserved electricity demand, the slack variables $\UnmetStorPlus$ and $\UnmetStorMinus$ for the storage level are added to the SPs, which are penalized by additional costs within the objective function (see Eq. (\ref{eq:objSPTS})). Therefore, the storage level is now calculated as
\begin{align}
\label{eq:storlevelTS}
    \StorLevel & = \StorLevelPrev + \Charge \etaStorIn - \frac{\Discharge}{\etaStorOut} \\ & \nonumber -\StorLoss + \UnmetStorPlus - \UnmetStorMinus, \\ & \nonumber \forall \text{ } \scenIndex \in \scenAll,\timeIndex \in \timeAll,\regionIndex \in \regionAll,\techStorage \in \techStorAll.
\end{align}
This ensures that the SPs remain feasible and only optimality cuts are added to the MP.

\subsection{Other enhancement strategies}
\label{subsec:acceleration}
\noindent As stated before, the classic BD algorithm may need a lot of iterations until it converges. Therefore, a variety of enhancement strategies have been proposed in the literature to improve the performance \cite{lumbreras2016solve,rahmaniani2017benders,goke2024stabilized}. Besides the adjustment of the decomposition strategy, we add further enhancement strategies to our algorithm, which are described in more detail in the following sections.

\subsubsection{MPI and GMI}
\label{subsubsec:MPI}
The SPs within the BD can be solved independently of each other, offering a high parallelization potential. Usually, when BD is calculated with shared memory, only one single node on the HPC system is used, limiting the parallelization potential. By implementing MPI (message passing interface) within the BD, one MPI process can be defined for each SP, which can then be solved on distributed memory, utilizing several computational nodes. This has the benefit that more resources can be used in parallel, allowing for a faster calculation process.
Furthermore, the model generation time can be quite time-consuming, especially within the large SPs. And since the model generation needs to be repeated in each iteration, this can lead to high time consumption for the same repeating process. This can be avoided by keeping the model open after each iteration and only updating the new information from the MP, similar to a sensitivity analysis where the model is kept in memory. Within GAMS, so-called “model instances” (GMI) can be used for this purpose \cite{Gmi}.

\subsubsection{Multi-cuts}
\label{subsubsec:multi}
In the classical BD algorithm, a single cut is generated in each iteration. For this, the dual variables of each SP are summed up (taking their respective probabilities into account) to a single cut. However, the dual variables of each SP can be considered in a separate cut, resulting in the multi-cut formulation. This approach has the benefit that more cuts are generated in each iteration, adding more information to the MP and resulting in a faster convergence of the algorithm \cite{jacobson2024computationally}. At the same time, the size of the MP increases more rapidly. Typically, this results in one cut per stochastic scenario. However, as described in Section \ref{subsec:timesplitting}, the multi-cut formulation can be extended when TSBD is applied. A SP is formulated for each scenario and time block combination, multiplying the number of cuts that can be added to the MP in each iteration by the number of time blocks.

\subsubsection{Bundle method}
\label{subsubsec:BM}
While it can be proven that the classical BD is able to find the optimal solution, this process can be very time-consuming, needing a lot of iterations until convergence. This issue is also described by Göke et al.\, where bundle methods are recommended to bundle the solution searching process within a trusted area (surrounding an initial starting solution and later close to the best solution found so far) \cite{goke2024stabilized}. They show that this can reduce the number of iterations needed until convergence significantly. As a starting point, one scenario of the stochastic problem can be picked, and the deterministic model can be solved. The bundle methods only allow the MP to search for new solutions within a limited radius. For this, a specific radius around the previous best solution can be picked, or the distance to the previous best solution can be penalized by additional costs within the objective function. In our case, we received the best performance when using the level bundle method \cite{frangioni2020standard}. The lower bound needs to stay below the level parameter
\begin{align}
\label{eq:lbmconstraint}
    \lbm \geq \Zexp + \sum_{\scenIndex,\timeBlock} \Thetatb
\end{align}
while minimizing the distance to the stability center, i.e.\ the capacities $\CapacityConvLBest$, $\CapacityStorLBest$ and $\CapacityTransLBest$ from the current best solution in iteration $\iterBest$, replacing Eq. (\ref{eq:objMPTS}) when solving the stabilized MP. The level parameter $\lbm$ is calculated as a weighted average between the lower and upper bound using the level weighting parameter $\beta$ \cite{goke2024stabilized}. The level bundle method performed best in our case when adjusted to minimize the distance of the summed capacities over all regions
\begin{align}
\label{eq:lbmobj}
     \min & \sum_{\techConverter} (\sum_{\regionIndex}(\CapacityConv - \CapacityConvLBest))^2 \\ & \nonumber + \sum_{\techStorage} (\sum_{\regionIndex}(\CapacityStor - \CapacityStorLBest))^2 \\ & \nonumber + \sum_{\techTransport} (\sum_{\regionIndex}(\CapacityTrans - \CapacityTransLBest))^2.
\end{align}
The stabilized MP becomes quadratic and, therefore, more complex to solve. To receive the actual lower bound of the objective value, another unstabilized linear MP is solved.

\subsubsection{Inactive cuts}
\label{subsubsec:inactivecut}
While considering cuts not only for each stochastic scenario but also for each time block, the number of cuts added per iteration to the MP increases considerably, as does the size and complexity of the MP. This can result in the MP becoming more time-consuming to solve than the SP after a certain number of iterations. Simultaneously, cuts generated in earlier iterations can become irrelevant for later iterations \cite{rahmaniani2017benders}. Therefore, cuts that were not binding for $\phi$ iterations in the stabilized and unstabilized MP can be deactivated for the next iterations \cite{goke2024stabilized}. The parameter $\delt$ stores the number of the iteration when the cut was created or lastly binding. The information, if a cut is active, is stored in the parameter $\act$, which is set to one after a cut is generated and to zero if a cut is deactivated. The cuts are considered not binding if their dual variable $\pibdcut$ is below a certain threshold $\epsactive$. This decreases the size of the MP again and its solving time. The convergence of the BD is not affected since relevant cuts that have been deactivated can be regenerated in later iterations. However, it is key to find the best fitting number of iterations when to deactivate a cut to avoid the need to regenerate too many cuts while still reducing the size and solving time of the MP.

Algorithm \ref{alg:benders_short} shows the TSBD using the regionally-summed level bundle method and inactive cuts presented in this section.

\begin{algorithm}[H]
\caption{Temporally split Benders decomposition using further enhancement methods.}
\label{alg:benders_short}
\begin{algorithmic}
\STATE \textbf{Input}: $\beta$, $\phi$, $\epsactive$, $\epsconv$
\STATE \textbf{Initialize}: $\iterBest \gets 1$, $\delt \leftarrow 0$, $\act \leftarrow 0$, $\lbm \leftarrow \inf$, $\Zlow \leftarrow \text{-} \inf$, $\Zup \leftarrow \inf, \Zupglob \leftarrow \inf$
\STATE Solve deterministic CEP to receive \textbf{starting point}
\STATE \textbf{fix} $\CapacityConv$, $\CapacityStor$, $\CapacityTrans$ and $\storagefixedVar$ in stab. MP in first iteration $\iterIndex=1$
\FOR{$\iterIndex \in \{1,...,\iterAll\}$}
  \STATE solve \textbf{stab. MP}
  \WHILE{stab. MP infeasible}
    \STATE $\Zlow \leftarrow \lbm$
    \STATE $\lbm \leftarrow \beta \Zlow + (1-\beta) \Zupglob$
    \STATE solve \textbf{stab. MP}
  \ENDWHILE
  \FOR{$\iterIndex' \in \{1, ...,\iterIndex-1\},\scenIndex \in \scenAll,\timeBlock \in \timeBlockAll$}
    \IF{$\pibdcutother > \epsactive$}
      \STATE $\deltother \leftarrow \iterIndex$
    \ENDIF
  \ENDFOR
  \STATE send $\CapacityConv$, $\CapacityStor$, $\CapacityTrans$ and $\storagefixedVar$ to SP
  \STATE solve \textbf{unstab. MP}
  \FOR{$\iterIndex' \in \{1, ...,\iterIndex-1\},\scenIndex \in \scenAll,\timeBlock \in \timeBlockAll$}
    \IF{$\pibdcutother > \epsactive$}
      \STATE $\deltother \leftarrow \iterIndex$
    \ENDIF
    \IF{$\iterIndex - \deltother > \phi$}
      \STATE $\actother \leftarrow 0$
    \ENDIF
  \ENDFOR
  \STATE fix $\CapacityConv,\CapacityStor,\CapacityTrans$ and $\storagefixedVar$ and solve each \textbf{SP} in parallel
  \STATE get $\piconvtb$, $\pistortb$, $\pitransaltb$, $\pitransagtb$, $\pistorfix$, $\pistorfixprev$, $\Zsptb$ and send to stab. and unstab. MP
  \STATE $\Zlow \leftarrow \Zmpunstab$
  \STATE $\Zup \leftarrow \Zexpstab + \sum_{\scenIndex,\timeBlock} \prob \Zsptb$
  \IF{$\Zup<\Zupglob$}
    \STATE $\Zupglob \leftarrow \Zup$
  \ENDIF
  \IF{$\Zup-\Zlow < \Zup_{k-1}-\Zlow_{k-1}$}
    \STATE $\iterBest \leftarrow \iterIndex$
  \ENDIF
  \STATE $\lbm \leftarrow \beta \Zlow + (1-\beta) \Zupglob$
  \STATE $\act \leftarrow 1$, $\delt \leftarrow \iterIndex$
  \IF{$1 - \Zlow/\Zupglob < \epsconv$}
    \STATE \textbf{exit for}
  \ENDIF
\ENDFOR
\end{algorithmic}
\end{algorithm}

\section{Case Study}
\label{sec:model}
\noindent The considered model is based on \cite{cao2018incorporating} and focuses on the power sector of Germany. One year is optimized with hourly resolution using a green-field optimization approach. The original dataset consists of 465 German nodes representing transformer substations. The imports from and exports to Germany's neighboring countries are considered with historical time series. However, the model can be spatially aggregated to any user-defined size, which facilitates the analysis of different levels of complexity. We assume that the power sector is fully decarbonized. Therefore, only renewable energies, biomass- and hydrogen-fueled power plants can be expanded. Open- and combined-cycle gas turbine power plants can be operated using green hydrogen, which can be imported for a price of 120 €/MWh$_{H2}$ \cite{lopion2020modellgestutzte}. Temporal and spatial balancing of supply and demand can be realized using pumped hydro, lithium-ion battery storage, and power transmission lines, respectively. Different uncertainties can be considered within the model \cite{sasanpour2022quantifying}. In this study, we consider weather uncertainties as random variables within the stochastic scenarios $\scenIndex$. For this purpose, we take seven historical weather years (2006 – 2012) with equal probability into account. Table \ref{tab:instance_size} lists the number of regions considered in the different model instances and the size of the DEQ taking the weather uncertainties into account. The largest instance consists of 30 aggregated German nodes and 9 nodes representing Germany's neighboring countries.

\begin{table}[t!]
	\caption{Model instances and size of deterministic equivalent with weather uncertainties.}
	\label{tab:instance_size}
    \small
	\begin{tabular}{rrrr}
		\hline
        Regions & Constraints [mio] & Variables [mio] & Non-zeros [mio] \\
        \hline
        \hline
        4 & 9.86 & 10.13 & 35.07\\
        13 & 17.49 & 19.00 & 62.68\\
        19 & 34.30 & 35.85 & 122.24\\
        39 & 87.04 & 88.98 & 309.14\\
        \hline
	\end{tabular}
\end{table}

\section{Results}
\label{sec:results}
\noindent The optimization with BD is performed on the HPC system CARO \cite{Caro} with 1276 standard nodes, each with two AMD~EPYC~7702 processors and 256~GB DDR4~memory and connected via a 100~GBit/s Infiniband network. The DEQ of the larger instances can not be solved on the standard nodes since the model runs out of memory. Therefore, we solve the DEQs on the 20 big memory nodes with 1024~GB DDR4 memory each. The DEQ, starting point, MP, and the SPs are solved with GAMS 48.2 and CPLEX 22.1 using the barrier method. The unstabilized MP is solved with dual simplex. For our calculations applying the (TS)BD algorithm, we use a parametrization of $\beta=0.1$, $\phi=50$, and $\epsactive=10^{-6}$ and a barrier convergence tolerance of $10^{-6}$. The convergence tolerance for the BD and the DEQs is set to $\epsconv=10^{-3}$.

\begin{table}[t!]
    \centering
	\caption{Optimized capacities of each model instance in GW.}
	\label{tab:capacities}
    \small
	\begin{tabular}{lrrrr}
		\hline
        Technology & 4r & 13r & 19r & 39r \\
        \hline
        \hline
        Biomass-fueled power plants & 9.5 & 8.4 & 10.6 & 11.1 \\
        H$_2$-fueled power plants & 54.3 & 54.0 & 57.7 & 58.9 \\
        Hydro run-of-river & 6.4 & 6.4 & 6.5 & 6.6 \\
        Photovoltaic & 228.0 & 227.8 & 237.7 & 241.9 \\
        Wind onshore & 69.1 & 68.7 & 67.3 & 62.9 \\
        Wind offshore & 0.2 & 0.3 & 12.3 & 20.1 \\
        Transmission grid & 283.0 & 344.0 & 510.8 & 843.6 \\
        Lithium-ion battery & 38.6 & 41.7 & 88.4 & 103.4 \\
        Pumped-hydro storage & 8.9 & 8.9 & 8.9 & 9.0 \\
        \hline
	\end{tabular}
\end{table}

The optimized capacities for the different model instances are listed in Table \ref{tab:capacities}. A higher spatial resolution results in higher total installed capacities. Especially the power grid is expanded to a significantly higher extend but also wind offshore becomes considerably more attractive.
Before comparing the performance of BD with and without temporal splitting, we analyze the impact of parallelizing the algorithm using MPI and restarting the SPs without regenerating them in each iteration using GMI (see Section \ref{subsubsec:MPI}). For this comparison, we solve the small- (4r) and medium-sized (13r) models with BD without applying temporal splitting.

\begin{figure}[!t]
    \centering
    \includegraphics[width=78mm]{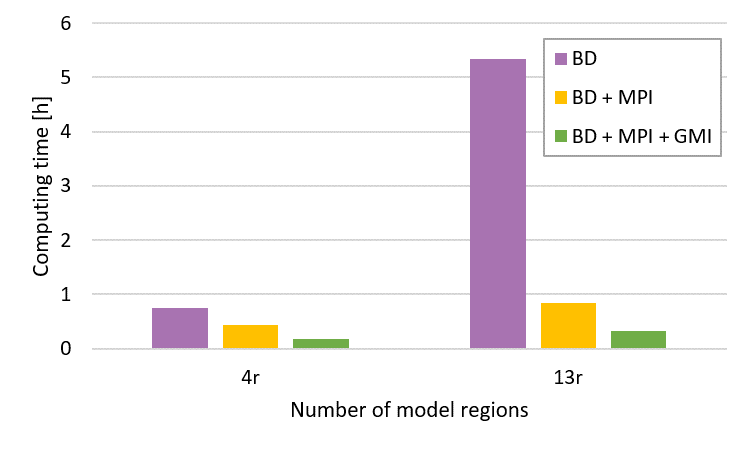}
    \caption{Comparison of computing time between different model sizes and different BD configurations, without applying temporal splitting.}
    \label{fig:Comparison}
\end{figure}

Fig. \ref{fig:Comparison} shows the computing time of the different BD configurations. The computing time refers to the total time, including model generation and solving. The "BD" configuration solves the models using shared memory. The 4r model is solved in 0.75~h, the 13r model in over 5h. The "BD + MPI" configuration makes use of the distributed memory of the HPC system and solves the problems in parallel on several computing nodes. This reduces the computing time by 42\% for the 4r model. The larger 13r model profits even more from the parallelization, where the computing time is reduced by 84\%. When using the "BD + MPI + GMI" configuration, we additionally leave each SP open, only updating the new capacities from the MP in each iteration. This additional feature further reduces the computing time by around 60\%, since the SPs only need to be generated once in the first iteration. Therefore, the 13r model can be solved with time savings of 94\% in total compared to the "BD" configuration. The 19r and 39r models can not be solved within 24h with the "BD" configuration, emphasizing the importance of parallelization.

\begin{figure}[!t]
    \centering
    \includegraphics[width=78mm]{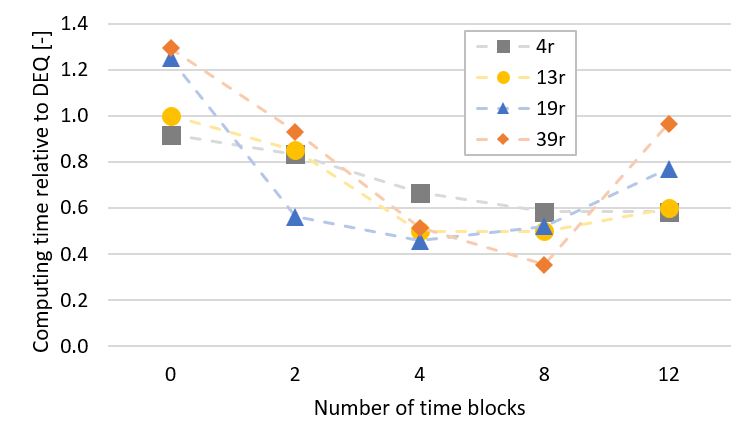}
    \caption{Computing time relative to computing time of DEQ.}
    \label{fig:Regions}
\end{figure}

Next, we compare the computing time of BD with and without temporal splitting. For this, we use the "BD + MPI + GMI" configuration due to the significant time reduction shown in Fig. \ref{fig:Comparison}.

In Fig. \ref{fig:Regions}, the computing times of BD with and without temporal splitting are compared to the computing time of the DEQ. Solving the 4r and 13r model with BD results in similar computing times as the DEQ, even without applying temporal decomposition (0 time blocks). The larger the model without temporal splitting, the higher the computing time compared to the DEQ. The 19r and 39r models are solved 20-30\% slower than the DEQ without temporal splitting. The computing time of all models decreases when TSBD is applied. While the small 4r model benefits the least, still a reduction in the computing time of up to 40\% can be achieved when decomposing the model into 8 or 12 time blocks. A decomposition into more time blocks results in lower computing times for the 4r model. The larger the model, the higher the relative time savings that can be achieved when applying the temporal splitting. The computing time of the largest 39r model is reduced by more than 60\% when applying a temporal splitting into 8 time blocks. For the mid- and large-sized models, the computing time increases again if 12 instead of 8 time blocks are selected. However, the larger models are more negatively impacted if the model is decomposed into too many time blocks. The computing time of the 39r model more than doubles if 12 instead of 8 time blocks are applied. Therefore, we take a closer look at the solving times of the 39r model.

\begin{figure}[!t]
    \centering
    \includegraphics[width=78mm]{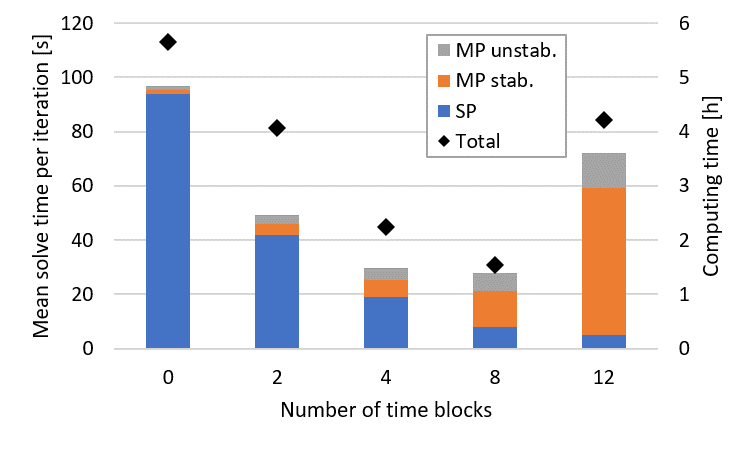}
    \caption{Mean solving time per iteration and total computing time for 39r model.}
    \label{fig:39r}
\end{figure}

Fig. \ref{fig:39r} shows the mean solving time per iteration and the total computing time for the 39r model, with and without temporal splitting. If no temporal splitting is applied (0 time blocks), the mean solving time for the first SP is much larger than the solving time of the stabilized and unstabilized MP. The mean SP solving time per iteration decreases significantly when temporal splitting is applied, and the time saving is higher with more selected time blocks. The mean MP solving time per iteration increases only slightly until 4 time blocks are applied. If a temporal splitting into 8 time blocks is applied, the stabilized MP is the most time-consuming component within the algorithm. Nevertheless, the total solving time per iteration decreases compared to applying 4 time blocks. Therefore, the computing time is also the lowest with 8 time blocks. With 12 time blocks selected, the mean solving time of the quadratic stabilized MP increases substantially, resulting in a more than doubled computing time for the calculation.

\begin{figure}[!t]
    \centering
    \includegraphics[width=78mm]{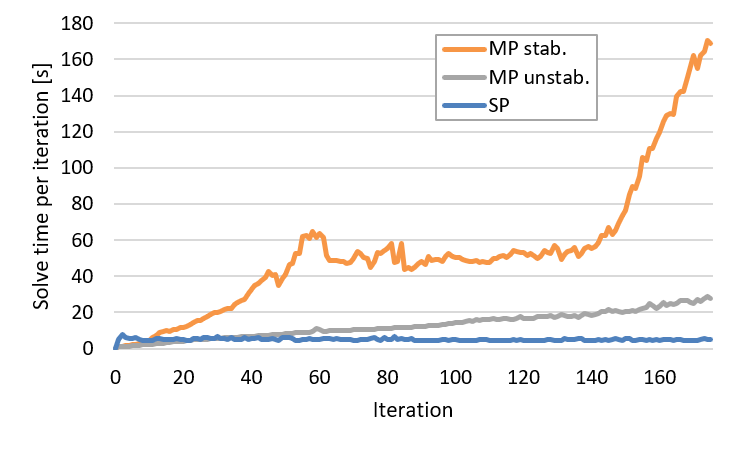}
    \caption{Solving time per iteration for 39r model and 12 time blocks applied.}
    \label{fig:12TB}
\end{figure}

Fig. \ref{fig:12TB} depicts the solving time of the first SP, the stabilized and unstabilized MP in each iteration for the 39r model and 12 time blocks applied. The solving time of the SP remains almost constant throughout the iterations, while the solving time of the unstabilized MP increases slightly and linearly. The solving time of the stabilized MP increases much faster than the unstabilized MP, with a high increase until iteration 50. Here, the first cuts can be deactivated and removed if they are not binding. The solving time of the stabilized MP remains almost constant until iteration 140, however, it fluctuates more. Afterwards, the solving time increases linearly again, but with a steeper ascent. A high share of cuts before iteration 90 are not binding and can be removed, decreasing the size of the MP again. However, cuts after iteration 90 become more relevant and are therefore not to be removed. The size of the MP increases again, making the quadratic problem more complex and time-consuming to solve.

\begin{figure}[!t]
    \centering
    \includegraphics[width=78mm]{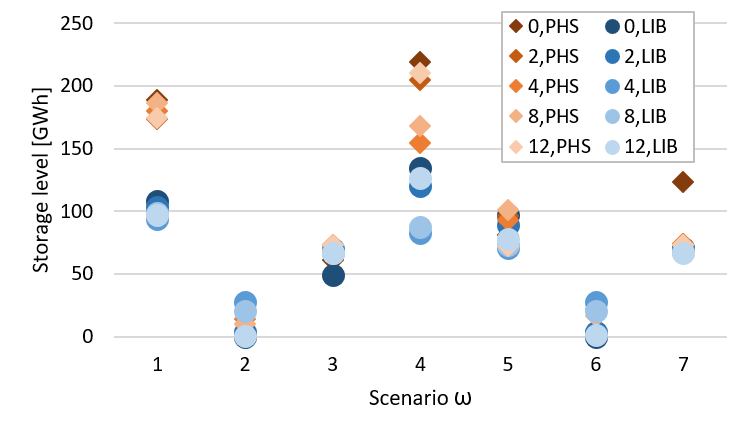}
    \caption{Storage level of 39r model in last time step, summed regionally for each scenario and different numbers of time blocks.}
    \label{fig:Storlevel}
\end{figure}

To analyze the effectiveness of the BDTS method and its optimization of the storage level in the last time step of each time block, we compare the storage level with and without temporal decomposition. Fig. \ref{fig:Storlevel} shows the storage level in the last time step for the 39r model and each scenario $\scenIndex$, representing different weather years. The storage levels of the pumped hydro storage (PHS) and the lithium-ion battery storage (LIB) are summed regionally. The analysis reveals that the storage level of PHS with BDTS, representing medium-term storage technologies, approximates the storage level without temporal decomposition. Despite the short-term storage duration of LIB, the varying storage levels of the seven stochastic scenarios can be accurately approached, irrespective of the number of time blocks selected. The larger deviations in the storage level observed in scenario 4 relate to the considered convergence tolerance of $\epsconv=10^{-3}$.

\section{Conclusion and Outlook} \label{sec:conclusions}
\noindent This paper introduces a novel method for parallelizing stochastic energy system models for capacity expansion planning with high temporal resolution. Through the implementation of a temporally split Benders decomposition, we manage to reduce the size of the subproblems representing the hourly energy system dispatch of the stochastic scenarios. In doing so, we enable a parallelization of these subproblems, which are usually much larger and more time-consuming to solve compared to the master problem. By optimizing the storage level of the last time step of each time block within the master problem, an optimal operation of long-term storage technologies can be achieved despite the temporal splitting. Our case study reveals that this approach can reduce model solution times by 60\% and lessen memory requirements compared to solving the deterministic equivalent.
To reduce the computing time even further, we combine the temporally split Benders decomposition with other enhancement strategies, such as bundle methods, extended multi-cuts, removal of inactive cuts, and solve the problems in parallel with distributed memory on high-performance computers using MPI. The application of distributed memory leads to further computing time savings of over 80\%. The four analyzed use-cases perform best when the problem is decomposed into 8 time blocks. Our results indicate that further increasing the number of time blocks is not favorable, as the solving time of the master problem increases significantly in later iterations due to the high number of added optimality cuts per iteration. A limitation is given by the comparatively slow convergence of the currently used bundle method in the first iterations, as it adds many unbinding cuts to the master problem. From this follows that future improvements of the bundle method could lead to a convergence in fewer iterations, decreasing the size of the master problem significantly, and improving the performance of the algorithm.

\appendices
\small{
\section*{CRediT author statement}
\noindent \textbf{Shima Sasanpour}: Methodology, Conceptualization, Investigation, Formal analysis, Writing - Original draft, Visualization.
\textbf{Manuel Wetzel}: Methodology, Formal analysis, Writing - Reviewing and Editing.
\textbf{Karl-Kiên Cao}: Data curation, Formal analysis, Writing - Reviewing and Editing.
\textbf{Hans Christian Gils}: Formal analysis, Writing - Reviewing and Editing, Funding acquisition.
\textbf{Andrés Ramos}: Supervision, Formal analysis, Writing - Reviewing and Editing.

\section*{Acknowledgment}
\noindent The research for this paper was performed within the projects 'UNSEEN' and 'ARTESIS' supported by the German Federal Ministry for Economic Affairs and Energy under grant numbers 03EI1004A and 03EI1067A. 

The authors gratefully acknowledge the scientific support and HPC resources provided by the German Aerospace Center (DLR). The HPC system CARO is partially funded by "Ministry of Science and Culture of Lower Saxony" and "Federal Ministry for Economic Affairs and Climate Action".

\ifCLASSOPTIONcaptionsoff
  \newpage
\fi

\bibliographystyle{IEEEtran}
\bibliography{references}

\end{document}